%%%%%%%%%%%%%%%%%%%%%%%%%%%%%%%%%%%%%%%%%%%%%%%%%%%%%%%%%%%%%
%%  Ramanujan's formula for the logarithmic derivative of 
%%                the gamma function           
%%
%%               author David Bradley 
%%
%%%%%%%%%%%%%%%%%%%%%%%%%%%%%%%%%%%%%%%%%%%%%%%%%%%%%%%%%%%%%

\input amstex
\documentstyle{amsppt}
\magnification=\magstephalf
\loadmsam
\loadmsbm

\TagsOnRight
\redefine\i{\infty}
\redefine\il{\int_0^{\i}}
\define\s{\sum_{k=1}^{\i}}

%\magnification=1200
%\NoBlackBoxes
%\NoRunningHeads

%\redefine\Re{\operatorname{Re}}
\define\a{\alpha}
\redefine\b{\beta}
\define\g{\gamma}
\define\z{\zeta}
\define\({\left(}
\define\){\right)}
\define\[{\left]}
\define\]{\right]}

\topmatter

\title Ramanujan's formula for the logarithmic derivative of 
the gamma function 
\endtitle

\rightheadtext{Ramanujan's formula for gamma'/gamma}

\author David Bradley \endauthor

\address Department of Mathematics, University of Illinois at
Urbana-Champaign \endaddress

\email dbradley\@math.uiuc.edu \endemail

\abstract We prove a remarkable formula of Ramanujan for the 
logarithmic 
derivative of the gamma function, which converges more rapidly
than classical expansions, and which is stated without proof in
the notebooks \cite5.  The formula has a number of very interesting
consequences which we derive, including an elegant hyperbolic summation,
Ramanujan's formula for the Riemann zeta function evaluated at the
odd positive integers, and new formulae for Euler's constant, $\g$.
\endabstract


%\footnote""{Presented at the 4th conference of the Canadian Number Theory
%Association, Halifax, 1994.}

\thanks Presented at the 4th conference of the Canadian Number Theory
Association, Halifax, 1994 \endthanks

\endtopmatter

\document
\interdisplaylinepenalty=500
%\baselineskip=24pt plus 2pt

%\noindent
Ramanujan was evidently fond of series expansions and representations of
special functions.  The unorganized material in the second and third
notebooks contain many interesting formulae involving the gamma
function, Bessel functions, and hypergeometric functions, to name a few
examples.  The following remarkable formula for the logarithmic
derivative of the gamma function has some very interesting consequences,
and should be contrasted with the well-known series representation 
$$ \psi(x+1) := \dfrac{\Gamma^{'}}{\Gamma}(x+1) = -\g +
\sum_{n=1}^{\i}\(\dfrac{1}{n}-\dfrac{1}{x+n}\),$$
[1, p.259, 6.3.16].
\medskip
\proclaim{Theorem}
For all $ x>0, $
$$ \align
\psi(x+1)=&\dfrac{\pi}{3}\log x+\dfrac{1}{2x}-\dfrac{1}{4\pi x^2}+
\dfrac{\pi\cot(\pi x)}{e^{2\pi x}-1} 
+\dfrac{\pi\log|2\sin(\pi x)|}{2\sinh^2(\pi x)} \\
&+\s\dfrac{2 k}{(e^{2\pi k}-1)(k^2-x^2)} 
-\dfrac{\pi}{2}\s\dfrac{\log|k^4-x^4|}{\sinh^2(\pi k)} \\
&-2\pi\s e^{-2\pi kx}\(k^2\sum_{n=1}^\i\dfrac{\sin(2\pi nx)}{k^2+n^2}-
k^3\sum_{n=1}^\i\dfrac{\cos(2\pi nx)}{n(k^2+n^2)}\).\tag1
\endalign
$$
\endproclaim
\medskip

%\noindent
The theorem is recorded as formula 2, p\. 280 in Ramanujan's
notebooks \cite5.  As per usual, Ramanujan gives no proof or explanation.
In the sequel, we give a proof, and offer a plausible explanation of
how Ramanujan may have arrived at formula (1).  But first, we examine the
formula more closely, with a view to gaining a better understanding
of the formula and its many interesting consequences, some of which 
may not be immediately apparent.

\medskip
The formula takes some time to absorb, but even an initial glance
reveals that the convergence is more rapid than the classical
representation.  The appearance of the logarithm in a summand is very
striking, if not mysterious, and at first may arouse some suspicion,
since $\psi(x+1)$ is well-behaved for $x>0$, whereas on the right we
have terms which contribute singularities at the positive
integers.  However, a closer examination reveals that the singularities
arising from the various terms actually cancel each other out.  For
example, from the partial fraction decomposition
$$ \pi\cot(\pi x) = \dfrac{1}{x} + \s\dfrac{2x}{x^2-k^2}$$
it can be seen that the pole contribution from 
$$ \dfrac{\pi\cot(\pi x)}{e^{2\pi x}-1} $$
at $x=m \quad (m\in \bold Z^{+})$ is cancelled by the $m$th term in the
sum
$$ \s\dfrac{2 k}{(e^{2\pi k}-1)(k^2-x^2)} $$
leaving only a harmless (finite) quantity.  This is a routine exercise.
Similarly, using the fact that 
$$ |\sin(\pi x)| \sim \pi |x-m| \qquad (x\rightarrow m) $$
one can show that the logarithmic singularity contributed by 
$ \dfrac{\pi\log|2\sin(\pi x)|}{2\sinh^2(\pi x)} $
at $x=m \quad (m\in \bold Z^{+})$ is cancelled by the $m$th term in the 
sum 
$$ \dfrac{\pi}{2}\s\dfrac{\log|k^4-x^4|}{\sinh^2(\pi k)} $$
Carrying out the details of the aforementioned limit calculations provides
the following formula for Euler's constant, $\g$.

\medskip
\proclaim{Corollary} Let $m$ be a positive integer.  Then 
$$ \align
-\g+\sum_{j=1}^m\dfrac{1}{j}=&\dfrac{\pi}{3}\log m+\dfrac{1}{2m}-
\dfrac{1}{4\pi m^2}
+\sum \Sb k=1 \\ k\neq m\endSb ^{\i}\dfrac{2 k}{(e^{2\pi k}-1)(k^2-m^2)} \\
&+\dfrac{\pi}{2}\(\dfrac{\log\pi -\log(2m^3)-1}{\sinh^2(\pi m)}\)
-\dfrac{\pi}{2}\sum \Sb k=1 \\ k\neq m \endSb ^{\i}\dfrac{\log|k^4-m^4|}
{\sinh^2(\pi k)} \\
&+\dfrac{1}{2 m\(e^{2\pi m}-1\)} 
+2\pi\s k e^{-2\pi k m} \sum_{n=1}^\i\dfrac{k^2}{n(k^2+n^2)}.\tag2
\endalign
$$
\endproclaim
\medskip

\noindent
Note that for the inner sum in (2) we have [1, p.259, 6.3.17]
$$ \sum_{n=1}^\i\dfrac{k^2}{n(k^2+n^2)} = \g + \Re \psi(1+i k).$$
We have used MAPLE
to check (2) with various values of $m$.  It is interesting to
note that with $m=1$, taking a mere five terms from each of the 
series indexed 
by $k$ is sufficient to give $1-\g$ accurate to thirteen decimal
places: $1-\g = 0.4227843350984\dots$.
\medskip
\noindent{\bf Proof of Corollary.} Let $x \rightarrow m$ in (1). 
It suffices to show that   
$$\lim_{x\to m} \(\dfrac{\pi\cot(\pi x)}{e^{2\pi x}-1}+
\dfrac{2m}{(e^{2\pi m}-1)(m^2-x^2)}\) = 
\dfrac{1}{2m(e^{2\pi m}-1)} - \dfrac{\pi}{2 \sinh^2(\pi m)}$$
and 
$$\lim_{x\to m}\(\dfrac{\pi\log|2\sin(\pi x)|}{2 \sinh^2(\pi x)}-
\dfrac{\pi \log|m^4-x^4|}{2\sinh^2(\pi m)}\) =
\dfrac{\pi\(\log \pi - \log(2m^3)\)}{2 \sinh^2(\pi m)}.$$
We have
$$\align \pi\cot(\pi x) =& \dfrac{\pi \cos(\pi(x-m))}{\sin(\pi(x-m))}
= \dfrac{\pi(1+O(x-m)^2)}{\pi(x-m)(1+O(x-m)^2)}\\
=& \dfrac{1}{x-m} + O(x-m),\endalign$$
$$\dfrac{1}{e^{2\pi x}-1} = \dfrac{1}{e^{2\pi m}-1} - \dfrac{2\pi
e^{2\pi m}}{(e^{2\pi m}-1)^2}(x-m) + O(x-m)^2,$$
and
$$\align \dfrac{2m}{(e^{2\pi m}-1)(m^2-x^2)} =& \dfrac{1}{(e^{2\pi
m}-1)(m-x)}\dfrac{1}{1+\frac{x-m}{2m}} \\
=& \dfrac{1}{(e^{2\pi
m}-1)(m-x)}\(1-\dfrac{x-m}{2m}+O(x-m)^2\) \\
=& \dfrac{1}{(e^{2\pi m}-1)(m-x)} + \dfrac{1}{2m(e^{2\pi
m}-1)}+O(x-m).\endalign$$
Thus,
$$\dfrac{\pi\cot(\pi x)}{e^{2\pi x}-1}+
\dfrac{2m}{(e^{2\pi m}-1)(m^2-x^2)} = 
\dfrac{1}{2m(e^{2\pi m}-1)} - \dfrac{2\pi e^{2\pi m}}{(e^{2\pi m}-1)^2}
+O(x-m),$$
which proves the former.  For the latter, we have
$$\align L :=& 
\lim_{x\to m}\(\dfrac{\log|2\sin(\pi x)|}{\sinh^2(\pi x)}-
\dfrac{\log|m^4-x^4|}{\sinh^2(\pi m)}\)\\
 =&
\dfrac{\log 2 - \log(4m^3)}{\sinh^2(\pi m)} + \lim_{x\to m}
\(\dfrac{\log|\sin(\pi x)|}{\sinh^2(\pi x)} - \dfrac{\log|m-x|} 
{\sinh^2(\pi m)}\) \\
=& -\dfrac{\log(2m^3)}{\sinh^2(\pi m)} + \lim_{x\to m}
\(\dfrac{\log\left|\pi(m-x)+O(m-x)^2\right|}{\sinh^2(\pi x)}-
\dfrac{\log|m-x|}{\sinh^2(\pi m)}\)\\
=& \dfrac{\log \pi - \log(2m^3)}{\sinh^2(\pi m)}+\lim_{x\to m}
\(\dfrac{\log|m-x|+O(m-x)}{\sinh^2(\pi x)}-
\dfrac{\log|m-x|}{\sinh^2(\pi m)}\)\\
=& \dfrac{\log \pi - \log(2m^3)}{\sinh^2(\pi m)},\endalign$$
as required.
\qed

\medskip
%\noindent
Instead of letting $x$ tend to a positive integer in (1), one can of
course obtain a formula for $\g$ by directly substituting 
$x = m - 1/2 \quad (m\in \bold Z^{+})$ in (1).  Jonathan Borwein 
has used this approach to
compute $\g$, choosing large $m$ to balance the error terms.  However,
he reports that Brent's method is computationally superior.  Part of the
difficulty stems from the need for full precision calculation of the
logarithms, and of $e^{\pi k}$.  If there is an efficient way of doing
these calculations, it is possible that an algorithm based on
Ramanujan's formula could beat other methods of computing $\g$. 
[J. Borwein, personal communication.]

\medskip 
%\noindent
Another very interesting, and perhaps unexpected consequence of (1) is
obtained by letting $x \rightarrow +\i$.  Since this limiting case is an
essential ingredient in the proof of (1), we shall temporarily content 
ourselves with heuristic reasoning, deferring the rigorous details of
the argument to the sequel.  As $x \rightarrow +\i$, most terms drop
out, and since $\psi(x+1) \sim \log x$, what remains is
$$ \log x \sim \dfrac{\pi}{3}\log x -
\dfrac{\pi}{2}\s\dfrac{\log|k^4-x^4|}{\sinh^2(\pi k)}
\qquad (x \rightarrow +\i).\tag3$$
For large $x$ and small $k$, $\log|k^4-x^4| \sim 4\log x$, whereas when
k is large, $\sinh^2(\pi k)$ is very large, so that (3) becomes
$$ \log x \sim \dfrac{\pi}{3}\log x - 2\pi\log
x\s\dfrac{1}{\sinh^2(\pi k)}\qquad (x \rightarrow
+\i).\tag4$$
For (4) to be true, the coefficients of $\log x$ must agree i.e. we
must have 
$$ \s\dfrac{1}{\sinh^2(\pi k)} = \dfrac{1}{6} -
\dfrac{1}{2\pi},\tag5$$
a remarkable result, with an interesting history in and of itself.
See [2, Prop\. 2.26].  That (5) is a special case of (1) makes 
Ramanujan's formula all the more intriguing.

\medskip 
%\noindent
An additional consequence of (1) is Ramanujan's remarkable formula 
for the Riemann
zeta function, evaluated at the odd positive integers.

\medskip
\proclaim{Corollary} Let $N$ be a positive integer.  Then
$$
\multline 2N\z(2N+1) + 4N\s \dfrac{k^{-2N-1}}{e^{2\pi k}-1} +
\(1+(-1)^N\)\pi\s \dfrac{k^{-2N}}{\sinh^2(\pi k)}\\
= (2\pi)^{2N+1}\sum_{j=0}^{N+1}(-1)^{j+1}(2j-1)\dfrac{B_{2j}}{(2j)!}
\dfrac{B_{2N+2-2j}}{(2N+2-2j)!}.
\endmultline\tag6
$$
\endproclaim
\noindent Here, $B_0=1$, $B_1=-1/2$, $B_2 = 1/6$, etc. are the Bernoulli
numbers, defined by 
$$
\dfrac{z}{e^z-1} = \s\dfrac{B_k}{k!}z^k, \qquad (|z|<2\pi).\tag7
$$
Since the Bernoulli numbers are all rational, (6) shows that $\z(2N+1)$ is
a rational multiple of $\pi^{2N+1}$ plus either one or two rapidly
convergent series, the number depending on the parity of $N$.
\medskip
\noindent{\bf Proof of Corollary.} The main idea stems from the
observation that $\psi(x+1)$ is essentially the ordinary generating
function for the sequence $\z(2), \z(3),\dots$.  More precisely, 
$$
\psi(x+1) = -\g + \sum_{n=1}^\i x^n (-1)^{n+1} \z(n+1), \qquad (|x|<1).
\tag8
$$
Thus, we are tempted to expand both sides of (1) into power series in
$x$ and then equate coefficients of $x^{2N}$.  However, since $\log x$
has no Maclaurin series, we first differentiate both sides of (1) with
respect to $x$.  Since the machinations involved in this differentiation
arise in our proof of (1), we shall defer the details to the sequel and
content ourselves for the moment with merely reporting that after
simplifying the resulting expression, we obtain
$$\align  
\psi^{'}(x+1) =& \dfrac{\pi}{3x} -\dfrac{1}{2x^2} + \dfrac{1}{2\pi x^3} 
-\dfrac{\pi^2\csc^2(\pi x)}{e^{2\pi x}-1} \\
&+\s\dfrac{4 k x}{(e^{2\pi k}-1)(k^2-x^2)^2}
+\s\dfrac{2\pi x^3}{\sinh^2(\pi k)(k^4-x^4)}.\tag9
\endalign
$$
We now expand all terms of (9) into power series in $x$.  We have
$$\align
\s\dfrac{4 k x}{(e^{2\pi k}-1)(k^2-x^2)^2} =& \s\dfrac{4x}{k^3\(e^{2\pi
k}-1\)}\sum_{n=1}^\i n\(\dfrac{x}{k}\)^{2n-2}\\
=& 4 \sum_{n=1}^\i nx^{2n-1} \s\dfrac{k^{-2n-1}}{e^{2\pi k}-1},\tag10 
\endalign
$$
$$\align
\s\dfrac{2\pi x^3}{\sinh^2(\pi k)(k^4-x^4)} =& \s\dfrac{2\pi
x^3}{k^4 \sinh^2(\pi k)}\sum_{n=0}^\i \(\dfrac{x}{k}\)^{4n} \\
=& 2\pi\sum_{n=1}^\i x^{4n-1}\s\dfrac{k^{-4n}}{\sinh^2(\pi
k)},\tag11
\endalign
$$
and from (8), 
$$
\psi^{'}(x+1) = \sum_{n=1}^\i (-x)^{n-1} n\z(n+1).\tag12
$$
Finally, using (7) and the derivative of the well-known power series for
$\pi \cot(\pi x)$, we have
$$\align
\dfrac{\pi^2\csc^2(\pi x)}{e^{2\pi x}-1} =& \dfrac{1}{2\pi
x^3}\sum_{n=0}^\i (2\pi x)^n\dfrac{B_n}{n!}\sum_{k=0}^\i
(-1)^{k+1}(2k-1)x^{2k}(2\pi)^{2k}\dfrac{B_{2k}}{(2k)!}\\
=& \dfrac{1}{2\pi x^3}\sum_{n=0}^\i (2\pi x)^{2n} \sum_{j=0}^n
(-1)^{j+1} (2j-1) \dfrac{B_{2j}}{(2j)!}\dfrac{B_{2n-2j}}{(2n-2j)!}\\
&+\dfrac{1}{2\pi x^3}\sum_{n=0}^\i (2\pi x)^{2n+1} \sum_{j=0}^n
(-1)^{j+1}(2j-1)\dfrac{B_{2j}}{(2j)!}\dfrac{B_{2n+1-2j}}{(2n+1-2j)!}\\
=& \dfrac{1}{2\pi x^3} +\dfrac{\pi}{3x} +\dfrac{1}{2\pi
x^3}\sum_{n=2}^\i (2\pi x)^{2n} \sum_{j=0}^n (-1)^{j+1} (2j-1)
\dfrac{B_{2j}}{(2j)!}\dfrac{B_{2n-2j}}{(2n-2j)!}\\
&- \dfrac{1}{2x^2} 
+\dfrac{1}{2\pi x^3}\sum_{n=1}^\i (2\pi x)^{2n+1} (-1)^{n+1} (2n-1)
\dfrac{B_{2n}}{(2n)!}\dfrac{B_1}{1!},
\tag13
\endalign
$$
since $B_{2n+1-2j} = 0$ for $j=0,1,\dots,n-1$.
Substituting (10) through (13) into (9) yields 
$$
\align
\sum_{n=1}^\i (-x)^{n-1} n\z(n+1) =& 4\sum_{n=1}^\i
nx^{2n-1}\s\dfrac{k^{-2n-1}}{e^{2\pi k}-1} +2\pi\sum_{n=1}^\i 
x^{4n-1}\s\dfrac{k^{-4n}}{\sinh^2(\pi k)}\\
&-\dfrac{1}{2\pi x^3}\sum_{n=2}^\i (2\pi x)^{2n} \sum_{j=0}^n (-1)^{j+1} 
(2j-1) \dfrac{B_{2j}}{(2j)!}\dfrac{B_{2n-2j}}{(2n-2j)!}\\
&-\dfrac{1}{2\pi x^3}\sum_{n=1}^\i (2\pi x)^{2n+1} (-1)^{n+1} (2n-1)
\dfrac{B_{2n}}{(2n)!}\dfrac{B_1}{1!}.
\tag14
\endalign
$$
The corollary (6) now follows on equating coefficients of $x^{2N-1}$ in
(14).  
\qed

\medskip \noindent
{\bf Remark 1.}  Equating coefficients of $x^{2N}$ in (14) gives Euler's
celebrated formula $\z(2N) = (-1)^N(2\pi)^{2N}\dfrac{B_{2N}}{(2N)!}
B_1 = (-1)^{N+1} 2^{2N-1} \pi^{2N}\dfrac{B_{2N}}{(2N)!}$,
which is also valid when $N=0$, since $\z(0)=-1/2 = B_1.$
It should also be noted that (6)
holds for $N=0$ if we make use of the hyperbolic summation (5) 
and interpret the term $2N\z(2N+1)$ as
$\dsize \lim_{x\to 0} 2x\z(2x+1)=1$.
In fact, it can be shown that (6) holds for negative integers $N$ as
well.  For example, the case $N=-1$ in (6) gives
$$
\s\dfrac{k}{e^{2\pi k}-1} = \dfrac{1}{24}-\dfrac{1}{8\pi}.
$$
As a further example, let $m>1$ be odd, and let $N=-m$.  From (6) and
the well-known evaluation $\z(1-2m) = -\dfrac{B_{2m}}{2m}$, we obtain
$$
\s\dfrac{k^{2m-1}}{e^{2\pi k}-1} = \dfrac{B_{2m}}{4m}.
$$
Recalling that $\dsize\il\dfrac{x^{2m-1}}{e^{2\pi x}-1}\,dx =
\dfrac{B_{2m}}{4m}$, we derive a rather striking equality between a 
sum and an integral in which the summand and the integrand are 
formally identical:
$$
\s\dfrac{k^{2m-1}}{e^{2\pi k}-1} = \il\dfrac{x^{2m-1}}{e^{2\pi x}-1}\,dx,
\qquad (\text{odd } m > 1).
$$

\medskip \noindent
{\bf Remark 2.}  Let $\a,\b>0$ with $\a\b=\pi^2$, and let $N$ be a
positive integer.  One can derive the slightly more general formula
$$
\multline 2N\a^{-N}\(\z(2N+1) + 2\s \dfrac{k^{-2N-1}}{e^{2\a k}-1}\)\\
+ \a^{1-N}\s \dfrac{k^{-2N}}{\sinh^2(\a k)} 
- (-\b)^{1-N}\s \dfrac{k^{-2N}}{\sinh^2(\b k)}\\
= 2^{2N+1}\sum_{j=0}^{N+1}(-1)^{j+1}(2j-1)\a^{N+1-j}\b^j\dfrac{B_{2j}}{(2j)!}
\dfrac{B_{2N+2-2j}}{(2N+2-2j)!}
\endmultline\tag15
$$
in the same manner as our proof of (6) by expanding $\dfrac{\pi
\sqrt{\pi\a}\csc^2(z\sqrt{\pi\a})}{e^{2z\sqrt{\pi\b}}-1}$ into
partial fractions and then comparing coefficients of like powers of
$z$.  The author is currently investigating the feasibility of this
relatively elementary approach to deriving analogous Ramanujan-like
formulae for the Dirichlet $L$-functions via related partial fraction
expansions.  

\medskip 
%\noindent
Other proofs of (15) have been given via Mellin transforms, 
the Poisson summation formula, and elliptic modular transformations, to
name some examples of various methods employed.
See [2] and [3, p\. 276] for an extensive bibliography.  
Since Ramanujan was quite fond of partial
fraction expansions, it is possible that our derivation is 
close to Ramanujan's.  

\medskip 
%\noindent
We have seen (2) that specializing $x$ to positive integer values in
(1) produces formulae for Euler's constant.  However, it should be 
noted that the following formula for Euler's constant, valid for all
real $x>0$, is actually equivalent to (1), namely,
$$ \align
-\g =& \dfrac{\pi}{3}\log x +\dfrac{1}{4\pi x^2}
+\dfrac{\pi\log|2\sin(\pi x)|}{2\sinh^2(\pi x)} 
-\s\dfrac{x^2}{k(k^2+x^2)} \\
&+\s\dfrac{2 k}{(k^2+x^2)(e^{2\pi k}-1)}
-\dfrac{\pi}{2}\s\dfrac{\log|k^4-x^4|}{\sinh^2(\pi k)} \\
&-2\pi\s e^{-2\pi kx}\(k^2\sum_{n=1}^\i\dfrac{\sin(2\pi nx)}{k^2+n^2}-
k^3\sum_{n=1}^\i\dfrac{\cos(2\pi nx)}{n(k^2+n^2)}\).\tag17
\endalign
$$
\medskip
\noindent To see the equivalence of (1) and (17), consider the partial 
fraction expansion of $\dsize \frac{\pi
\cot(\pi x)}{e^{2 \pi x}-1}$ in the form
$$
\aligned
&\ \psi(x+1) + \g = \frac{1}{2x} - \frac{1}{2\pi x^2} + \frac{\pi 
\cot(\pi x)}{e^{2\pi x}-1}\\
&\ \hphantom{.....................} 
+\s{\frac{x^2}{k(k^2+x^2)}} 
+\s{\frac{4 k x^2}{(e^{2\pi k}-1)(k^4-x^4)}}
\endaligned\tag18
$$
(Entry 8, Chapter 30 [4, p\. 374]).  Subtracting (18) from (1), 
we obtain (17).   

\medskip 
%\noindent
If we use the classical formulae [1, p.259, 6.3.13, 6.3.17]
$$ \Im\psi(1+ix) = \dfrac{\pi}{2}\coth(\pi x) - \dfrac{1}{2x} \qquad 
\text{and} \qquad \Re\psi(1+ix) = -\g + \s\dfrac{x^2}{k(k^2+x^2)}$$
and substitute the latter in (17), after cancelling $-\g$ from
both sides there results a formula for
$\psi(1+ix)$, valid for $x>0$, analogous to formula (1) for $\psi(1+x)$.

\medskip 
%\noindent
The evaluation (17) was originally proved by Bruce Berndt, Jonathan
Borwein, and Will Galway.   They showed that the derivatives of both sides
of (17) 
vanish, and that the right side of (17) tends to $-\g$ as $x$ tends
to infinity.  Since (1) and (17) are equivalent, this argument provides
a verification of Ramanujan's formula.  
The question remains as to how Ramanujan arrived at (1) in the first
place.  We may never know how Ramanujan arrived at some of his
results, but in the case of (1), we are able to give a plausible
answer to our question.

\medskip 
%\noindent
Ramanujan was in the habit of expanding meromorphic functions into 
their series of partial fractions, and many such examples can be found
in the unorganized pages of the second notebook.  However, the following
partial fraction expansion (Entry 3 of Chapter 30 [4, p\. 359])
has a particular relevance to (1).
$$\align  
\dfrac{\pi^2\csc^2(\pi x)}{e^{2\pi x}-1} =& \dfrac{\pi}{3x} 
-\dfrac{1}{2x^2} + \dfrac{1}{2\pi x^3} -\psi^{'}(x+1) \\
&+\s\dfrac{4 k x}{(e^{2\pi k}-1)(x^2-k^2)^2}
-\s\dfrac{2\pi x^3}{\sinh^2(\pi k)(x^4-k^4)}.\tag19
\endalign
$$
Observe that (19) is simply a rearrangement of (9).
One is tempted to speculate that Ramanujan, observing the presence of
the derivative $\psi^{'}(x+1),$ was motivated to integrate both sides of 
(19) to obtain a formula for $\psi(x+1)$.  Carrying out the
anti-differentiation reveals a happy correspondence between the terms
that arise and several of the terms appearing on the right of (1).
To indicate just one example, we see the mysterious 
$\dfrac{\pi}{2}\dsize\s\dfrac{\log|k^4-x^4|}{\sinh^2(\pi k)}$
appearing as the anti-derivative of $\dsize\s\dfrac{2\pi x^3}
{\sinh^2(\pi k)(x^4-k^4)}.$  However, a few terms on the right of
(1) remain unaccounted for, of which the double series is the most
glaringly obvious.  We shall see that the double series
arises fairly naturally in the course of proving (1).  Accordingly, we
now focus on completing this final task.

\medskip
\noindent{\bf Proof of Theorem.} 
By the product rule for differentiation, we have
$$
\frac{d}{dx}\dfrac{\pi\cot(\pi x)}{e^{2\pi x}-1} = 
- \dfrac{\pi^2 \csc^2(\pi x)}{e^{2\pi x}-1} 
- \dfrac{\pi^2 \cot(\pi x)}{2 \sinh^2(\pi x)}
$$
and
$$
\frac{d}{dx}\dfrac{\pi\log|2\sin(\pi x)|}{2\sinh^2(\pi x)} =
\dfrac{\pi^2 \cot(\pi x)}{2 \sinh^2(\pi x)}
+\dfrac{\pi}{2}\log|2\sin(\pi x)|\frac{d}{dx}\text{csch}^2(\pi x).
$$
Hence, we can rewrite (19) in the form
$$\align  
\psi^{'}(x+1) =& \dfrac{\pi}{3x} -\dfrac{1}{2x^2} + \dfrac{1}{2\pi x^3} 
+ \frac{d}{dx}\dfrac{\pi\cot(\pi x)}{e^{2\pi x}-1} 
+ \frac{d}{dx}\dfrac{\pi\log|2\sin(\pi x)|}{2\sinh^2(\pi x)}\\
&+ \s\dfrac{4 k x}{(e^{2\pi k}-1)(x^2-k^2)^2}
- \s\dfrac{2\pi x^3}{\sinh^2(\pi k)(x^4-k^4)}\\
&- \dfrac{\pi}{2}\log|2\sin(\pi x)|\frac{d}{dx}\text{csch}^2(\pi x)
.\tag20
\endalign
$$
But,
$$\align
- \dfrac{\pi}{2}\frac{d}{dx}\text{csch}^2(\pi x)
=& -\frac{d}{dx}\dfrac{2\pi e^{2\pi x}}{\(e^{2\pi x}-1\)^2}
= \(\frac{d}{dx}\)^2\dfrac{1}{e^{2\pi x}-1}
= \(\frac{d}{dx}\)^2\s e^{-2\pi kx}\\
=& \s (-2\pi k)^2 e^{-2\pi kx}.
\endalign
$$
Using this in (20) and integrating both sides with respect to $x$, we
obtain
$$ \align
\psi(x+1)=&\dfrac{\pi}{3}\log x+\dfrac{1}{2x}-\dfrac{1}{4\pi x^2}+
\dfrac{\pi\cot(\pi x)}{e^{2\pi x}-1} 
+\dfrac{\pi\log|2\sin(\pi x)|}{2\sinh^2(\pi x)} \\
&+\s\dfrac{2 k}{(e^{2\pi k}-1)(k^2-x^2)} 
-\dfrac{\pi}{2}\s\dfrac{\log|k^4-x^4|}{\sinh^2(\pi k)} \\
&- \int_x^\i \log|2\sin(\pi v)|\s (-2\pi k)^2 e^{-2\pi kv}\,dv + C
,\tag21
\endalign
$$
where $C$ is a constant of integration to be determined.  We now let
$$
\align
S(x) :=& \int_x^\i \log|2\sin(\pi v)|\s (-2\pi k)^2 e^{-2\pi kv}\,dv\\
=& \il\log|2\sin(\pi(x+u))|\s (-2\pi k)^2e^{-2\pi k(x+u)}\,du. 
\tag22
\endalign
$$
By means of the Laplace transforms
$$ \il e^{-kt}\sin(nt)dt = \dfrac{n}{k^2+n^2} \qquad \text{and} \qquad
\il e^{-kt}\cos(nt)dt = \dfrac{k}{k^2+n^2}, $$
and the well--known Fourier series
$$ \s\dfrac{\cos(2\pi ky)}{k}= -\log|2\sin(\pi y)|, \qquad y \text{ real}, $$
we find from (22) that
$$
\align
S(x) =& (2\pi)^2\s k^2e^{-2\pi kx}
         \il e^{-2\pi ku}\log|2\sin(\pi(x+u))|\,du,\\
=& -(2\pi)^2\s k^2e^{-2\pi kx}\il e^{-2\pi ku}\sum_{n=1}^\i
         \dfrac{\cos(2\pi nx + 2\pi nu)}{n}\,du\\
=& -2\pi\s k^2e^{-2\pi kx}\il e^{-kt}\sum_{n=1}^\i
         \dfrac{\cos(2\pi nx + nt)}{n}\,dt\\
=& 2\pi\s k^2e^{-2\pi kx}\il e^{-kt}\sum_{n=1}^\i
         \dfrac{\sin(2\pi nx)\sin(nt)- \cos(2\pi nx)\cos(nt)}{n}\,dt \\
=&  2\pi\s k^2e^{-2\pi kx}\(\sum_{n=1}^\i\dfrac{\sin(2\pi nx)}{n} 
         \il e^{-kt}\sin(nt)\,dt\right.\\  &\left. \qquad \qquad 
         \qquad \qquad \qquad -\sum_{n=1}^\i\dfrac{\cos(2\pi nx)}{n}
         \il e^{-kt}\cos(nt)\,dt\)\\
=& 2\pi\s e^{-2\pi nx}\(k^2\sum_{n=1}^\i\dfrac{\sin(2\pi nx)}{k^2+n^2}-
         k^3\sum_{n=1}^\i\dfrac{\cos(2\pi nx)}{n(k^2+n^2)}\),
\tag23
\endalign
$$
where the inversions in the order of summation and integration are justified
by absolute convergence.  Substituting (23) for the definition of $S(x)$
in (21), we have
$$ \align
\psi(x+1)=&C+\dfrac{\pi}{3}\log x+\dfrac{1}{2x}-\dfrac{1}{4\pi x^2}+
\dfrac{\pi\cot(\pi x)}{e^{2\pi x}-1} 
+\dfrac{\pi\log|2\sin(\pi x)|}{2\sinh^2(\pi x)} \\
&+\s\dfrac{2 k}{(e^{2\pi k}-1)(k^2-x^2)} 
-\dfrac{\pi}{2}\s\dfrac{\log|k^4-x^4|}{\sinh^2(\pi k)} \\
&- 2\pi\s e^{-2\pi nx}\(k^2\sum_{n=1}^\i\dfrac{\sin(2\pi nx)}{k^2+n^2}-
         k^3\sum_{n=1}^\i\dfrac{\cos(2\pi nx)}{n(k^2+n^2)}\).\tag24
\endalign
$$
By letting $x\rightarrow +\i$ in (24), we can evaluate the integration
constant, $C$.  In fact, $C$ must equal zero, as a comparison of (24) and
(1) reveals.  Since most terms in (24) vanish as $x\to+\i$, it
suffices to show that 
$$\left. \lim_{x\to\i} \( \psi(x+1) -\dfrac{\pi}{3}\log x + 
\dfrac{\pi}{2}\s\dfrac{\log|x^4-k^4|}{\sinh^2(\pi k)}\) = 0. \right.
\tag25
$$
Set $ x = N+1/2, $ where
$ N $ is a positive integer, and write
$$\align 
\s\dfrac{\log|x^4-k^4|}{\sinh^2(\pi k)}&=
\(\sum_{k\leq \sqrt{x}}+\sum_{\sqrt{x}<k<x}+\sum_{k>x}\)
\dfrac{\log|x^4-k^4|}{\sinh^2(\pi k)} \\
=&\sum_{k\leq \sqrt{x}}\dfrac{4\log x}{\sinh^2(\pi k)} 
+ \sum_{k\leq \sqrt{x}}\dfrac{\log(1-k^4/x^4)}{\sinh^2(\pi k)} \\
&+\(4x\log x\) O\(e^{-2\pi\sqrt{x}}\)+O\(e^{-3\pi x/2}\) \\
=&\(4\log x\)\(\s\dfrac{1}{\sinh^2(\pi k)}+O\(e^{-\pi\sqrt{x}}\)\) \\
&+O(x^{-2}) +O\(e^{-\pi\sqrt{x}}\) \\
=&\(4\log x\)\s\dfrac{1}{\sinh^2(\pi k)}+O(x^{-2}). 
\tag26
\endalign
$$
Second, by Stirling's formula,
$$ \psi(x+1) = \log x +O(1/x), \qquad (x\to+\i).\tag27 $$
Employing (26) and (27), we see that
$$
\align \psi(x+1) &-\dfrac{\pi}{3}\log x +
\dfrac{\pi}{2}\s\dfrac{\log|x^4-k^4|}{\sinh^2(\pi k)} \\
=&\(1-\dfrac{\pi}{3} + 2\pi\s\dfrac{1}{\sinh^2(\pi k)}\)\log x + 
O(1/x),\tag28
\endalign 
$$
as $x$ tends to $+\i$.  However (Berndt [2, Prop\. 2.26]),
$$
1-\dfrac{\pi}{3}+2\pi\s\dfrac{1}{\sinh^2(\pi k)} = 
1-\dfrac{\pi}{3}+2\pi\(\dfrac{1}{6}-\dfrac{1}{2\pi}\)=0.
\tag29 
$$
Using (29) in (28), we find that
$$ 
\psi(x+1)-\dfrac{\pi}{3}\log x +  
\dfrac{\pi}{2}\s\dfrac{\log|x^4-k^4|}{\sinh^2(\pi k)} = O(1/x),
\tag30 $$
as $ x $ tends to $ +\i. $  Thus, the limit (25) holds, and the proof of
the theorem (1) is complete. 
\qed

\medskip

\subhead Acknowledgement \endsubhead The author is grateful to Bruce
Berndt, Jonathan Borwein, and Will Galway for their helpful observations.


\bigskip

\Refs

%\ref\no1 \by Richard Askey and Mourad Ismail
%\paper A Combinatorial Sum
%\paperinfo MRC Technical Summary Report \#1557, July, 1975. 
%Unpublished
%\endref

\ref\no1 \by M. Abramowitz and I. Stegun
\book Handbook of Mathematical Functions
\publ Dover
\publaddr New York
\yr 1972
\endref

\ref\no2 \by Bruce C. Berndt
\paper Modular transformations and generalizations of several formulae
of Ramanujan
\jour Rocky Mt. J. Math. 
\vol 7
\yr1977
\pages 147-189
\endref

\ref\no3 \by Bruce C. Berndt
\book Ramanujan's Notebooks Part II
\publ Springer-Verlag
\publaddr New York
\yr1989
\endref

\ref\no4 \by Bruce C. Berndt
\book Ramanujan's Notebooks Part IV
\publ Springer-Verlag
\publaddr New York
\yr1994
\endref

\ref\no5 \by S. Ramanujan
\book Notebooks (2 Volumes)
\publ Tata Institute of Fundamental Research
\publaddr Bombay
\yr1957
\endref

\endRefs
\enddocument